\documentclass[12pt]{article}
\usepackage{latexsym,amsfonts,amsthm,amsmath,amscd,amssymb}
\usepackage[dvips]{graphicx}

\newtheorem{lemma}{Lemma}[section]
\newtheorem{thm}[lemma]{Theorem}
\newtheorem{rem}[lemma]{Remark}
\newtheorem{prop}[lemma]{Proposition}
\newtheorem{cor}[lemma]{Corollary}

\newcommand\matP{{\mathbb{P}}}
\newcommand\matR{{\mathbb{R}}}
\newcommand\matQ{{\mathbb{Q}}}
\newcommand\matN{{\mathbb{N}}}
\newcommand\matZ{{\mathbb{Z}}}

\renewcommand{\hbar}{{\overline{h}}}

\newcommand{\crn}{{\rm crn}}
\newfont{\Got}{eufm10 scaled 1200}

\newcommand{\mycap} [1] {\caption{\footnotesize{#1}}}

\newcommand{\vol}{\mathop{\rm vol}\nolimits}

\newcommand{\Tor}{\mathop{\rm Tor}\nolimits}

\newcommand\calS{{\mathcal S}}

\newcommand{\deffont}[1]{{\textit{\textbf{#1}}}}

\begin{document}

\title{Complexity of links in $3$-manifolds}

\author{Ekaterina~\textsc{Pervova}\and Carlo~\textsc{Petronio}}

\maketitle

\begin{abstract}
\noindent We introduce a complexity $c(X)\in\matN$ for pairs
$X=(M,L)$, where $M$ is a closed orientable $3$-manifold and
$L\subset M$ is a link. The definition employs simple spines, but
for well-behaved $X$'s we show that $c(X)$ equals the minimal number
of tetrahedra in a triangulation of $M$ containing $L$ in its
$1$-skeleton. Slightly adapting Matveev's recent theory of
\emph{roots} for graphs, we carefully analyze the behaviour of $c$
under connected sum away from and along the link. We show in
particular that $c$ is almost always additive, describing in detail
the circumstances under which it is not. To do so we introduce a
certain (0,2)-\emph{root} for a pair $X$, we show that it is
well-defined, and we prove that $X$ has the same complexity as its
$(0,2)$-root. We then consider, for links in the $3$-sphere, the relations
of $c$ with the crossing number
and with the hyperbolic volume of the exterior, establishing various
upper and lower bounds. We also
specialize our analysis to certain infinite families of links,
providing rather accurate asymptotic estimates.

\end{abstract}

\section*{Introduction}

When one wants to analyze a certain class of low-dimensional
topological objects, having in mind in particular the exploitation
of computers to enumerate and classify them, one is naturally faced
with the need of two tools as follows:
\begin{itemize}
\item[(1)] An efficient combinatorial encoding of the objects under investigation;
\item[(2)] For each object, a numerical measure of how complicated the object is.
\end{itemize}
For closed orientable $3$-manifolds, item (1) is realized taking
triangulations. These would also work for item (2), but it actually
turned out over the time that Matveev's theory of
complexity~\cite{MaCompl} defined through simple spines is highly
preferable, for the following reasons:
\begin{itemize}
\item Matveev's complexity $c(M)\in\matN$ is defined for each $M$, and, if $M$ is prime,
either $M\in\left\{S^3,\ S^2\times S^1,\ \matP^3,\ L(3,1)\right\}$,
or $c(M)$ is the minimal number of tetrahedra in a triangulation of
$M$;
\item $c$ is additive under connected sum;
\item While performing a computer enumeration of the prime $M$'s with
a given $c(M)$, the extra flexibility deriving from the definition of $c$
through spines provides very useful computational shortcuts.
\end{itemize}

Besides closed orientable $3$-manifolds, perhaps the next more natural objects one can
want to understand and enumerate within $3$-dimensional topology are
the pairs $X=(M,L)$, where $M$ is an arbitrary closed orientable
$3$-manifold and $L\subset M$ is a link. This paper is devoted to
developing a complexity theory for such pairs, mimicking Matveev's
one. The definition of the complexity $c(X)$ is actually a variation of that already given
in~\cite{PeOrbCompl} for $3$-orbifolds, and exploited in~\cite{HHMP}
for pairs $(M,G)$, with $G$ a trivalent graph. However, specializing
to links we will be able to establish more accurate results, and to
compare complexity with other known invariants.

One of the main tools we will employ is a reduction of Matveev's
recent beautiful \emph{root} theory~\cite{MaRoots}. This theory
applies to pairs $(M,G)$, where $G$ is again a graph, and its aim is
to provide a version for such pairs of the celebrated
Haken-Kneser-Milnor theorem, asserting that each closed $3$-manifold
can be uniquely expressed as a connected sum of prime ones. Roughly
speaking, a root of $(M,G)$ is the result of cutting $(M,G)$ as long
as possible along \emph{essential} spheres intersecting $G$
transversely in at most $3$ points. Matveev's main achievement
in~\cite{MaRoots} is then the proof that the root is essentially
unique.

Our main results are as follows:
\begin{itemize}

\item We will derive from~\cite{PeOrbCompl} the fact that for many pairs $X=(M,L)$
the complexity
$c(X)$ defined through simple spines is actually the minimal number
of tetrahedra in a triangulation of $M$ containing $L$ in its
$1$-skeleton; this is the case for instance if
$M=S^3$ and $L$ is a prime non-split link;

\item We will define a $(0,2)$-root as the result of cutting a pair $(M,L)$
as long as possible along essential spheres meeting $L$ in either
$0$ or $2$ points, and we will prove that such a root is
essentially unique;

\item We will prove that there are infinitely many distinct pairs $(M,L)$
in which every separating sphere meeting $L$ in either $0$ or $2$ points
is inessential, but containing non-separating spheres meeting $L$
in $2$ points; this shows in our opinion that there is
a non-negligible difference between root theory for pairs $(M,L)$ and
the Haken-Kneser-Milnor theorem for manifolds;

\item We will show that the complexity of $(M,L)$ is always equal to the sum of the complexities
of the components of its $(0,2)$-root;

\item We will prove that complexity is always additive under connected sum away
from the link, and it is additive also under connected sum along the link except in a situation
that we will carefully describe;

\item For links in $S^3$, we will compare the complexity with the crossing number $\crn$,
proving in particular for prime alternating knots that $c$ is bounded from above
by a linear function of $\crn$, and from below by a linear function
of $\log\left(\crn\right)$; we will also analyze two specific
infinite series of knots, for which we can give much sharper upper
and lower estimates;

\item We will compare the complexity of a link with the hyperbolic volume of
its exterior, showing that for certain families of links these two
quantities are rather closely related, while they cannot be in
general.

\end{itemize}

\bigskip

\noindent \textsc{Acknowledgements}: The first named author was
supported by the Marie Curie fellowship MIF1-CT-2006-038734. The
second named author is grateful to the CTQM in {\AA}rhus, the MFO in
Oberwolfach, and the University of Paris 7 for travel support.

\section{Links, spines, and complexity}

In this section we introduce the objects we will be dealing with throughout the
paper, and we define their complexity,
adapting the notion originally introduced in~\cite{PeOrbCompl} as a
variation of that given in~\cite{MaCompl}.

\paragraph{Spines of link-pairs}
We call \deffont{link-pair} a pair $X=\left(M,L\right)$, where
$M=:M\left(X\right)$ is a closed orientable 3-manifold and
$L=:L\left(X\right)$ is a (possibly empty) link in
$M\left(X\right)$. We will regard manifolds up to homeomorphism and
subsets of manifolds up to isotopy, without further mention. A
polyhedron $P$ (in the piecewise-linear sense~\cite{RoSa}) is called
\deffont{simple} if it is compact and the link of each point of $P$ embeds in
the complete graph with $4$ vertices. If the link of $p\in P$ is the
whole complete graph with $4$ vertices then $p$ is called a
\deffont{vertex} of $P$. A subpolyhderon $P$ of $M\left(X\right)$ is
called a \deffont{spine} of a link-pair $X$ if $P$ intersects
$L\left(X\right)$ transversely and the complement of $P$ in
$M\left(X\right)$ consists of balls either disjoint from
$L\left(X\right)$ or containing a single unknotted arc of
$L\left(X\right)$. The \deffont{complexity} $c\left(X\right)$ of $X$
is the minimal number of vertices of a spine of $X$.

\paragraph{Special spines, triangulations, and duality}
A polyhedron $P$ is called \deffont{quasi-special} if the link of
each point of $P$ is either a circle, or a $\theta$-graph, or the
complete graph with $4$ vertices. A quasi-special polyhedron $P$ is
of course simple. In addition, every point of $P$ either is a
vertex, or it belongs to a triple line, or it lies on a surface
region. This gives a stratification of $P$ into a $0$-, a $1$-, and
a $2$-dimensional set, and $P$ is called \deffont{special} if this
stratification is actually a cellularization.

A \deffont{triangulation} of  a pair $X$ is a realization of $M\left(X\right)$ as a gluing of a finite number of
tetrahedra along simplicial maps between the triangular faces, with the property that $L\left(X\right)$ is
contained in the union of the glued edges. The next result is easily established (see~\cite{PeOrbCompl}
for details):

\begin{prop}\label{duality:prop}
Associating to a triangulation of a link-pair $X$ the $2$-skeleton
of the cellularization of $M\left(X\right)$ dual to it one gets a
bijection between the set of triangulations of $X$ and the set of
special spines of $X$.
\end{prop}

\paragraph{\textit{i}-spheres and \textit{i}-irreducible pairs}
Given a pair $X$, we call \deffont{i-sphere} in $X$ a subset
of $M\left(X\right)$ homeomorphic to $S^2$ and meeting $L\left(X\right)$ transversely in $i$ points. In the sequel
we will be interested in the cases $i=0,1,2$ (and later $i=3$, when using the results of Matveev~\cite{MaRoots}).
For $i=0,1,2$ an $i$-sphere is \deffont{trivial} in $X$ if it bounds in $M\left(X\right)$ a ball
which is disjoint from $L\left(X\right)$ if $i=0$ and meets $L$ in a single unknotted arc when $i=2$ (so
a $1$-sphere is never trivial). The pair $X$ is called \deffont{i-irreducible} if it does
not contain non-trivial $i$-spheres.

\paragraph{Minimal spines and naturality of complexity}
A simple spine $P$ of a link-pair $X$ is called \deffont{minimal} if it has $c\left(X\right)$ vertices and no proper subset of $P$
is a spine of $X$. The next result is a direct consequence of~\cite[Theorem 2.6]{PeOrbCompl}:

\begin{thm}\label{natural:thm}
Let $X$ be a $\left(0,1,2\right)$-irreducible link-pair. Then:
\begin{itemize}
\item If $c\left(X\right)=0$ then $M\left(X\right)$ is either $S^3$, $L\left(3,1\right)$, or $\matP^3$, and $L\left(X\right)$
is either empty or the core of a Heegaard torus of
$M\left(X\right)$;
\item If $c\left(X\right)>0$ then each minimal spine of $X$ is special.
\end{itemize}
\end{thm}

\noindent
(The assumption of $3$-irreducibility in~\cite[Theorem 2.6]{PeOrbCompl}
is actually never used in the proof.)

\begin{cor}
With only $6$ exceptions (having complexity $0$), the complexity of
a $\left(0,1,2\right)$-irreducible link-pair $X$ is the minimal
number of tetrahedra in a triangulation of $X$.
\end{cor}

\section{Review of Matveev's root theory}
We will review in this section the theory developed in~\cite{MaRoots}, slightly changing
the terminology to adhere to that used above, and citing the results we will explicitly need.
This theory applies to \deffont{graph-pairs}, namely to pairs $X=\left(M,G\right)$ where $M=:M\left(X\right)$ is a closed
orientable $3$-manifold and $G=:G\left(X\right)$ is a unitrivalent graph embedded in $M\left(X\right)$.
Empty graphs and graphs with knot components are allowed, so link-pairs are also graph-pairs.
Arbitrary graphs, with vertices
of valence other than $1$ or $3$, may be accepted, but this would not add much to the theory.

\paragraph{Trivial \textit{i}-balls and \textit{i}-pairs, and (trivial) \textit{i}-spheres}
We will now slightly extend and adjust some of the definitions given above.
For $i=0,1,2,3$ we call \deffont{trivial i-ball} the pair
$\left(B^3,G_i\right)$ where $G_i$ is the cone over $i$ points of $\partial B^3$ with vertex
at the centre of $B^3$. Note that this is not strictly a graph-pair according
to the above definition, since it has boundary. We next define the \deffont{trivial i-pair}
to be the double of the trivial $i$-ball, mirrored in its boundary.

Given a graph-pair $X$ and $i\leqslant3$ we define an
\deffont{i-sphere} in $X$ to be a sphere embedded in
$M\left(X\right)$ that meets $G\left(X\right)$ transversely in $i$
points ---in particular, none of these points can be a vertex of
$G\left(X\right)$--- and we define $\Sigma$ to be a \deffont{trivial
i-sphere} in $X$ if it bounds a trivial $i$-ball in $X$. (Note
that if $G\left(X\right)$ is a link then any $1$-sphere is
non-trivial, coherently with above.)

\paragraph{Splittable and essential spheres}
Given a sphere $\Sigma$ in a graph-pair $X$ we call
\deffont{compression} of $\Sigma$ the operation of taking a disc $D$
in $M\left(X\right)$ such that $D\cap\Sigma=\partial D$ and $D\cap
G\left(X\right)=\emptyset$, and replacing $\Sigma$ by the two
boundary components of a regular neighbourhood of $\Sigma\cup D$
that do not cobound with $\Sigma$ a product $\Sigma\times[0,1]$. We
then say that an $i$-sphere $\Sigma$ is \deffont{splittable} if
there exists a compression of $\Sigma$ yielding an $i_1$-sphere and
an $i_2$-sphere with $i_1,i_2\geqslant 1$. Note that $i_1+i_2=i$, so
if $i\leqslant 1$ then an $i$-sphere is always unsplittable.

We say that an $i$-sphere in $X$ is \deffont{essential} if it is
unsplittable and non-trivial, and we say that $X$ is
\deffont{i-irreducible} if it does not contain essential
$i$-spheres. Note that a link-pair can be $2$-reducible as a
link-pair, as defined in the previous section, without being
$2$-reducible as a graph-pair as just defined. However the combined
definition of $\left(0,1,2\right)$-irreducibility, which is what we
will really need, is coherent.

\paragraph{Surgery}
If $\Sigma$ is an $i$-sphere in a graph-pair $X$ we call
\deffont{surgery on X along $\mathbf{\Sigma}$} the operation
of cutting $X$ open along $\Sigma$ and attaching two trivial
$i$-balls matching the intersections with the graph. Note that
$\Sigma$ may be separating or not. We will denote the result of such
a surgery by $X_\Sigma$.

The definition of surgery naturally extends to the case where instead of $\Sigma$
one takes a disjoint union $\calS$ of $i$-spheres (with varying $i$). Such
a disjoint union is called a \deffont{sphere-system}, and the result of the surgery
on $X$ along it is again denoted by $X_\calS$. The next fact will
be repeatedly used in the sequel:

\begin{prop}\label{system:for:subsequent:prop}
If a graph-pair $Y$ is the result of a finite number of subsequent
surgeries starting from some graph-pair $X$, then there exists a
sphere-system $\calS$ in $X$ such that $Y=X_\calS$.
\end{prop}

Note that the sphere-system $\calS$ given by this proposition is not
well-defined up to isotopy in $X$, but we will only use the fact that it exists.

\paragraph{(0,1,2,3)-roots}
The following was established in~\cite{MaRoots} using normal surfaces:

\begin{thm}\label{finite:surg:teo}
Given a graph-pair $X$, there exists $n\in\matN$ with the following
property: for every sequence $X=X_0,X_1,\ldots,X_k$ of graph-pairs
such that for each $j\geqslant 1$ the pair $X_j$ is obtained from
$X_{j-1}$ by surgery along an essential
$\left(0/1/2/3\right)$-sphere, one has $k\leqslant n$.
\end{thm}

We now define a graph-pair $Y$ to be a
\deffont{$\mathbf{(0,1,2,3)}$-root} of a pair $X$ if $Y$ is
$\left(0,1,2,3\right)$-irreducible and $Y=X_k$ for some sequence
$X=X_0,X_1,\ldots,X_k$ as in the previous statement. The statement
itself implies that $\left(0,1,2,3\right)$-roots exist for every
$X$. To any root $Y$ of $X$ we will always associate a sphere-system
$\calS$ such that $Y=X_\calS$, even if $\calS$ is actually not
uniquely defined up to isotopy.

\paragraph{Efficient (0,1,2,3)-roots}
We say that a sphere-system $\calS$ in $X$ is \deffont{$\mathbf{(0,1,2,3)}$-efficient} if
surgery on $X$ along $\calS$ gives a $\left(0,1,2,3\right)$-root of $X$ and for all $\Sigma\in\calS$
one has that $\Sigma$ is essential in $X_{\calS\setminus\Sigma}$.
A $\left(0,1,2,3\right)$-root of $X$ is in turn called \deffont{efficient}
if it is the result of a surgery along an efficient sphere-system.

As remarked in~\cite{MaRoots}, if one starts from a sphere-system
giving a $\left(0,1,2,3\right)$-root and one removes from it a
sphere violating the efficiency condition, then
the resulting sphere-system still gives a root of $X$.
This easily implies that efficient sphere-systems exist for every
$X$, and hence efficient $\left(0,1,2,3\right)$-roots also do.

\paragraph{Sliding moves and uniqueness of efficient (0,1,2,3)-roots}
Given two sphere-systems $\calS$ and $\calS'$, we say that $\calS'$
is obtained from $\calS$ by a \deffont{sliding} if the following
happens:
\begin{itemize}
\item There exist spheres $\Sigma\in\calS$ and $\Sigma'\in\calS'$ such that
$\calS\setminus\Sigma=\calS'\setminus\Sigma'$;
\item There exist an $i$-sphere $\Omega$ in $\calS\setminus\Sigma=\calS'\setminus\Sigma'$, with
$i\in\{0,2\}$, and a simple closed arc $\alpha$ in $X$ meeting
$\calS$ only at its ends one of which belongs to $\Sigma$ and the
other to $\Omega$, with $\alpha\cap
G\left(X\right)=\emptyset$ for $i=0$ and $\alpha\subset
G\left(X\right)$ for $i=2$;
\item $\Sigma'$ is the boundary component of the regular neighbourhood of
$\Sigma\cup\alpha\cup\Omega$ that does not cobound with $\Sigma$ or
$\Omega$ a product $\Sigma\times[0,1]$ or $\Omega\times[0,1]$.
\end{itemize}
It is very easy to see that under these assumptions one has $X_{\calS'}=X_\calS$.
This fact and the next result we state imply the main achievements of~\cite{MaRoots},
namely uniqueness of efficient $\left(0,1,2,3\right)$-roots and virtual uniqueness of arbitrary
$\left(0,1,2,3\right)$-roots:

\begin{thm}\label{sliding:thm}
Any two $\left(0,1,2,3\right)$-efficient systems of the same graph-pair are
related by a finite sequence of slidings.
\end{thm}

\begin{cor}\label{unique:ess:root:cor}
The efficient $\left(0,1,2,3\right)$-root is well-defined for each
graph-pair.
\end{cor}

\begin{cor}\label{quasi:unique:root:cor}
Any two $\left(0,1,2,3\right)$-roots of the same graph-pair are
obtained from each other by insertion and removal of trivial
graph-pairs.
\end{cor}

As a conclusion, we note that Corollary~\ref{unique:ess:root:cor}
can be viewed as an analogue for graph-pairs of the Haken-Kneser-Milnor
theorem, according to which every compact $3$-manifold can be uniquely expressed
as a connected sum of prime ones. The fact that to get a root one cuts also
along non-separating spheres however makes root theory somewhat weaker
than the decomposition into primes, even if stronger than the theory
developed in~\cite{PeOrbDeco}. The reasons are as follows:
\begin{itemize}
\item Two distinct pairs can have the same efficient root (take for
instance $N=M\#\left(S^2\times S^1\right)$, with empty links; then $N$ has the same efficient root as $M$);
\item A graph-pair $X$ may not be the result of a connected sum of the components
of its efficient root (which happens for the same $N$ as in the previous example);
\item There are infinitely many graph-pairs that are not $\left(0,1,2,3\right)$-irreducible such
that every separating $\left(0/1/2/3\right)$-sphere is trivial (examples were provided in~\cite{PeOrbDeco}).
\end{itemize}

\section{$\mathbf{\left(0,2\right)}$-reduction of Matveev's root theory}
We will extensively use in this section the terminology and results of the previous one.

\paragraph{(0,2)-roots}
We call \deffont{$\mathbf{(0,2)}$-root} of a link-pair
$X=\left(M,L\right)$ a $\left(0,2\right)$-irreduci\-ble pair obtained
from $X$ by subsequent surgery along essential
$\left(0/2\right)$-spheres. Theorem~\ref{finite:surg:teo} already
implies that every $X$ has $\left(0,2\right)$-roots, and again we
will use the fact that any such root can be obtained as $X_\calS$
for a system $\calS$ of $\left(0/2\right)$-spheres, even if $\calS$ is not
well-defined up to isotopy. We will say that such a
system $\calS$ is \deffont{efficient} if $\Sigma$ is essential in
$X_{\calS\setminus\Sigma}$ for each component $\Sigma$ of $\calS$,
and in this case we call $X_\calS$ an \deffont{efficient $\mathbf{(0,2)}$-root}.

\begin{prop}\label{eff:02syst:prop}
The $\left(0,2\right)$-efficient systems in a link-pair $X$ are precisely
those that one can obtain from a $\left(0,1,2,3\right)$-efficient system
of $X$ by discarding the $\left(1/3\right)$-spheres.
\end{prop}

\begin{proof}
Let $\calS$ be a $\left(0,2\right)$-efficient system in $X$. We first extend
$\calS$ to a system $\calS'$ giving a $\left(0,1,2,3\right)$-root of $X$, adding each time
a sphere that is essential after surgery along the previous ones, and
then we extract from $\calS'$ an efficient $\left(0,1,2,3\right)$-system $\calS''$.
We must show the following:
\begin{enumerate}
\item The only $\left(0/2\right)$-spheres contained in $\calS'$ are those of $\calS$;
\item No sphere from $\calS$ gets discarded when passing from $\calS'$ to $\calS''$.
\end{enumerate}
Item 1 is easy: if there exists a $\left(0/2\right)$-sphere
$\Sigma\in\calS'\setminus\calS$ then $\Sigma$ is essential
in $X_{\calS\cup\calS'''}$, for some $\calS'''\subset\calS'\setminus\calS$,
whence it is essential in $X_{\calS}$, which contradicts the fact that $X_\calS$ is
a $\left(0,2\right)$-root of $X$.

Turning to item 2, we argue by contradiction and let $\Sigma$ be the
first $\left(0/2\right)$-sphere discarded in passing from $\calS'$
to $\calS''$. Setting $Y=X_{\calS\setminus\Sigma}$ we have that
$\Sigma$ is essential in $Y$ but not in $Y_{\calS''''}$, where
$\calS''''$ is a system of $\left(1/3\right)$-spheres (those of
$\calS'\setminus\calS$ that have not been discarded yet). Suppose
first that $\Sigma$ is splittable in $Y_{\calS''''}$. This means
that there is a disc $D$ disjoint from $L\left(Y_{\calS''''}\right)$
and compressing $\Sigma$ into two $1$-discs. Moving $D$ away from
the traces of $\calS''''$ in $Y_{\calS''''}$ we can assume that $D$ lies
in $Y$ and is disjoint from $L(Y)$, so $\Sigma$ is
splittable in $Y$. This is a contradiction.

Therefore we can assume that $\Sigma$ is trivial in $Y_{\calS''''}$,
so it bounds a trivial $\left(0/2\right)$-ball there. However such a ball does
not contain any trivial $\left(1/3\right)$-ball, so it is disjoint from the traces
of ${\calS''''}$ in $Y_{\calS''''}$.
Therefore it can be assumed to lie in $Y$, which again gives a contradiction.

We have proved so far that efficient $\left(0,2\right)$-systems extend to
efficient $\left(0,1,2,3\right)$-system adding $\left(1/3\right)$-spheres only. Now we will prove that
given an efficient $\left(0,1,2,3\right)$-system $\calS$ the subsystem $\calS'$ consisting
of the $\left(0/2\right)$-spheres is an efficient $\left(0,2\right)$-system. Of course
for each $\Sigma$ in $\calS'$ we have that $\Sigma$ is essential in
$X_{\calS'\setminus\Sigma}$ because it is after the further surgery along
$\calS\setminus\calS'$. We are then left to show that $X_{\calS'}$
is a $\left(0,2\right)$-root of $X$, namely that it is $\left(0,2\right)$-irreducible. Suppose it is not,
and let $\Sigma$ be an essential $\left(0/2\right)$-sphere in $X_{\calS'}$. The proof
of item 2 above now shows that $\Sigma$ remains essential also in
$X_\calS$, because $\calS\setminus\calS'$ contains $\left(1/3\right)$-spheres only.
This gives a contradiction and the proof is complete.
\end{proof}

\paragraph{Uniqueness of efficient (0,2)-roots}
We now turn to the $(0,2)$-analogues of the main results of Matveev~\cite{MaRoots}:

\begin{prop}\label{exists:unique:eff02root:prop}
Any link-pair admits a unique efficient $\left(0,2\right)$-root.
\end{prop}

\begin{proof}
Proposition~\ref{eff:02syst:prop} and existence of $(0,1,2,3)$-efficient roots
readily imply existence. In addition, Proposition~\ref{eff:02syst:prop} together with
Theorem~\ref{sliding:thm} implies that any two efficient $\left(0,2\right)$-systems
are related by slidings. This is because a sliding always takes place along a $\left(0/2\right)$-sphere,
and does not change the type of the other sphere involved, hence a sliding
between two efficient $\left(0,1,2,3\right)$-systems refines to one on their subsystems
consisting of the $\left(0/2\right)$-spheres. The conclusion now follows from the remark
that sliding does not affect the result of a surgery.
\end{proof}

\begin{prop}\label{quasi:unique:02:root:prop}
Any two $\left(0,2\right)$-roots of one link-pair are obtained from each other
by insertion and removal of trivial link-pairs.
\end{prop}

\begin{proof}
Let $\calS_1$ and $\calS_2$ be $\left(0,2\right)$-systems such that
$X_{\calS_1}$ and $X_{\calS_2}$ are roots of $X$. Extract efficient
$\left(0,2\right)$-systems $\calS'_j\subset\calS_j$. Then
$X_{\calS'_1}=X_{\calS'_2}$. Now suppose $\calS'_j\setminus\calS_j$
contains an $i$-sphere $\Sigma$. Then $\Sigma$ is inessential in
$X_{\calS'_j}$, and it cannot be splittable otherwise it would be in
$X$. So it is trivial, which readily implies that
$X_{\Sigma'_j\cup\Sigma}$ is obtained from $X_{\Sigma'_j}$ by adding
a trivial $i$-pair. The conclusion follows by iterating this
argument.
\end{proof}

At the end of the previous section, in the context of graph-pairs,
we have listed three facts, making the point that root theory does not strictly
provide a version for graph-pairs of the Haken-Kneser-Milnor unique
decomposition theorem into primes for manifolds. The same facts stated
there also hold true for link-pairs, the first two of which
again by referring to the example of $N=M\#\left(S^2\times S^1\right)$.
For the third fact (existence of infinitely many non-irreducible primes)
the examples in~\cite{PeOrbDeco} do not work (they involve graphs with vertices); therefore we
introduce the pairs $X_n$ with $M(X_n)=S^2\times S^1$ as suggested in Fig.~\ref{Xn:fig},
and we prove the following:
    \begin{figure}
    \begin{center}
    \includegraphics[scale=0.8]{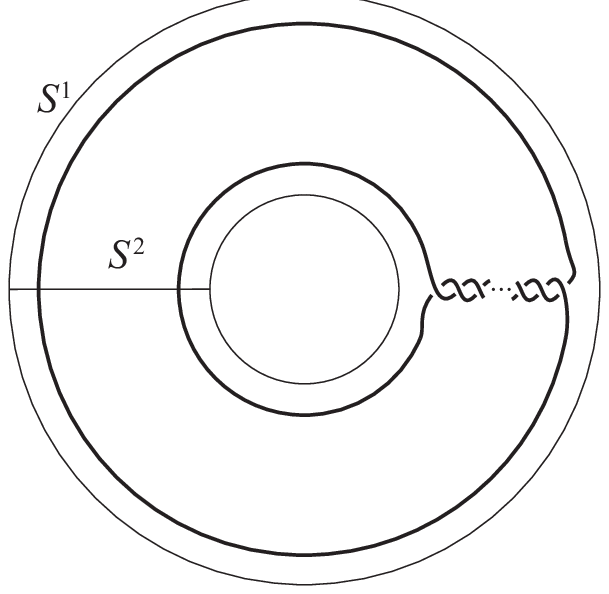}
    \mycap{The pair $X_n$, with $n\geqslant 1$ standing for the number of crossings in this picture.
    Note that $L(X_n)$ is always a knot}
    \label{Xn:fig}
    \end{center}
    \end{figure}

\begin{prop}\label{Xn:even:prop}
If $n\neq m$ then $X_n\neq X_m$. Moreover for even $n$
the pair $X_n$ is $\left(0,1\right)$-irreducible and not $2$-irreducible,
but every separating $2$-sphere in $X_n$ is trivial.
\end{prop}

\begin{proof}
The obvious $2$-sphere $\Sigma=S^2\times \{*\}$ is essential in $X_n$, and surgery along
it gives $(S^3,T_{2,n})$ where $T_{2,n}$ is the $(2,n)$ torus knot or link.
Since $T_{2,n}$ is a prime non-split link, $(S^3,T_{2,n})$ is $(0,2)$-irreducible,
so it is the efficient root of $X_n$, but the efficient root is unique by
Proposition~\ref{exists:unique:eff02root:prop}, and the first assertion follows.
Using Proposition~\ref{quasi:unique:02:root:prop} we also deduce that
\begin{itemize}
\item[$(\clubsuit)$] any $\left(0,2\right)$-root of $X_n$ contains precisely one non-trivial pair,
\end{itemize}
which we will need below.
We now consider the following property for a link-pair $X$:
\begin{itemize}
\item[$(*)$] $M(X)=S^3$ and $L(X)$ is a knot
\end{itemize}
and we make the obvious but crucial remark that
\begin{itemize}
\item[$(\diamondsuit)$] if $X$ satisfies $(*)$ then so does each component of a root of $X$.
\end{itemize}
For the rest of the proof we assume that $n$ is even, so $T_{2,n}$ has $2$ components, whence
\begin{itemize}
\item[$(\heartsuit)$] the only non-trivial component of a root of $X_n$ does \emph{not} satisfy $(*)$.
\end{itemize}

Let us now show that $X_n$ is $0$-irreducible. Suppose it is not, and let
$\Sigma$ be an essential $0$-sphere. Then any
$\left(0,2\right)$-root of $(X_n)_{\Sigma}$ is also a
$\left(0,2\right)$-root of $X_n$. If $\Sigma$ is non-separating then
$(X_n)_\Sigma$ satisfies $(*)$, so we get a contradiction by $(\diamondsuit)$
and $(\heartsuit)$. If $\Sigma$ is separating then $(X_n)_\Sigma$ is
the union of $S^2\times S^1$ and a pair satisfying $(*)$, but the
root of $S^2\times S^1$ is trivial, so again we get a contradiction
by $(\diamondsuit)$ and $(\heartsuit)$.

Of course $X_n$ cannot contain $1$-spheres, for otherwise $L(X_n)$ would
be non-trivial in $H^2(S^2\times S^1;\matZ)=H_1(S^2\times S^1;\matZ)$, whereas it is trivial.

We are left to show that in $X_n$ there is no non-trivial separating
$2$-sphere. Suppose that there is one, say $\Sigma$. Since
$X_n$ does not contain $1$-spheres, $\Sigma$ is essential, hence any
$\left(0,2\right)$-root of $(X_n)_{\Sigma}$ is also a
$\left(0,2\right)$-root of $X_n$. Now note that
$(X_n)_\Sigma=Y\sqcup Z$, where $M(Y)=S^2\times S^1$ and $M(Z)=S^3$,
and both $L(Y)$ and $L(Z)$ are knots. By $(\clubsuit)$, either $Y$
or $Z$ has a root consisting of trivial pairs only. If this happens
for $Z$ then $Z$ is itself trivial, because it satisfies $(*)$,
therefore $\Sigma$ is trivial in $X_n$, against our assumptions. We
conclude that a root of $Y$ consists of trivial pairs only, so the
only non-trivial pair from $(\clubsuit)$ is contained in a root of
$Z$. However $Z$ satisfies $(*)$ and once again we get a
contradiction from $(\diamondsuit)$ and
$(\heartsuit)$.
\end{proof}

The next result shows that one cannot remove the assumption that $n$ should be even
in the previous proposition:

\begin{prop}\label{Xn:odd:prop}
If $n\geqslant 3$ is odd then $X_n$ contains essential separating $2$-spheres.
\end{prop}

\begin{proof}
Let $\Sigma$ be the obvious non-separating $2$-sphere $S^2\times\{*\}$.
Since $n$ is odd, cutting $X_n$ along $\Sigma$ (without capping) we
get  $S^2\times [0,1]$ with two arcs each joining $\Sigma\times\{0\}$ to $\Sigma\times\{1\}$.
Take $\Sigma'$ to be the boundary of a regular neighbourhood of the union of $\Sigma$ and one
of these arcs. Then $\Sigma'$ is a separating $2$-sphere and
one easily sees that $(X_n)_{\Sigma'}$ is the union of $X_1$ and $(S^3,T_{2,n})$,
so $\Sigma'$ is essential.
\end{proof}

\section{Computation of complexity from roots,\\ and (restricted) additivity}
This section is devoted to the proofs of our main results about
complexity of link-pairs $X$ with $M(X)$ a general $3$-manifold
---starting from the next section we will mostly deal with the case $M(X)=S^3$.
Namely we will show that a link-pair has the same complexity as any
of its roots, and that complexity of link-pairs is additive under
connected sum away from the link and, provided the involved pairs do
not contain $1$-spheres, along the link. We will also show that when
there are $1$-spheres, additivity under connected sum along the link
does not hold. The notion of connected sum was already implicitly
referred to above, but we now define it formally.

\paragraph{\textit{i}-connected sum}
Let $X_1$ and $X_2$ be link-pairs. For $i\in\{0,2\}$ we define
a pair $X$ to be an \deffont{i-connected sum} of $X_1$ and $X_2$,
in symbols $X=X_1\#_i X_2$,
if $X$ is obtained as follows from $X_1$ and $X_2$:
\begin{itemize}
\item For $i=0$, remove from $M\left(X_j\right)$ a ball disjoint from $L\left(X_j\right)$, and glue the resulting
boundary $0$-spheres;
\item For $i=2$, remove from $M\left(X_j\right)$ a ball containing a single unknotted arc of $L\left(X_j\right)$,
and glue the resulting boundary $2$-spheres matching the intersections with the links.
\end{itemize}

\begin{rem}
\emph{If $M\left(X_1\right)$ and $M\left(X_2\right)$ are connected,
$X_1\#_0X_2$ can be performed in at most two ways,
and $X_1\#_2X_2$ in at most $4n_1n_2$ ways, where $n_j$ is the number
of components of $L\left(X_j\right)$.}
\end{rem}

\paragraph{Two-sided complexity estimates under surgery}
We will now examine what happens to complexity when doing surgery along an $i$-sphere,
showing that it can never diminish and that it stays the same when
the $i$-sphere satisfies a suitable assumption.

\begin{prop}\label{superadd:prop}
Let $\Sigma$ be a $\left(0/2\right)$-sphere in a link-pair $X$.
Then $c\left(X\right)\leqslant c\left(X_\Sigma\right)$.
\end{prop}

\begin{proof}
Let $P$ be spine of $X_\Sigma$ having $c\left(X_\Sigma\right)$ vertices.
We can enlarge $P$ a bit (without adding vertices) so that
$P\cap L\left(X_\Sigma\right)$ consists of surface points of $P$. Now note that the trace of
$\Sigma$ in $X_\Sigma$ consists of two spheres bounding trivial balls of
the appropriate type, and we can shrink these balls as much as we want.
Therefore we can assume they lie in $X_\Sigma\setminus P$. Adding bubbles to $P$
we can also assume that these balls:
\begin{itemize}
\item Lie in different components of $X_\Sigma\setminus P$
---this would be automatic if $\Sigma$ were separating in $X$ or if the traces of
$\Sigma$ in $X_\Sigma$ were incident to different components of
$L\left(X_\Sigma\right)$;
\item Are disjoint from $L\left(X_\Sigma\right)$ if $\Sigma$ is a $0$-sphere.
\end{itemize}
Then we get a spine of $X$ with $c\left(X_\Sigma\right)$ vertices by
adding to $P$ a tube (as in Fig.~9 of~\cite{PeOrbCompl}, for the
case where $\Sigma$ is a $2$-sphere), and the proposition is proved.
\end{proof}

\paragraph{Normal essential \textit{i}-spheres}
We recall here that a theory of normal 2-suborbifolds with respect to
handle decompositions of 3-orbifolds was developed in~\cite{PeOrbDeco}.
This theory applies \emph{verbatim} to link-pairs
and surfaces transverse to the link.
Given a simple spine $P$ of a pair $X$, if we triangulate $P$ (in the strict
piecewise-linear sense), we inflate each $j$-simplex of $P$ to a $j$-handle,
and we add 3-handles corresponding to the components of $M\left(X\right)\setminus P$,
we get a handle decomposition of $X$. We then say that a surface
transverse to $L\left(X\right)$ is in \deffont{normal position with respect to $\mathbf{P}$}
if it is with respect to a handle decomposition of $X$ induced by $P$.

\begin{prop}\label{subadd:prop}
Let $\Sigma$ be a $\left(0/2\right)$-sphere in a link-pair $X$
lying in normal position with respect to a spine $P$ of $X$
having $c\left(X\right)$ vertices. Then $c\left(X_\Sigma\right)\leqslant c\left(X\right)$.
\end{prop}

\begin{proof}
As in~\cite{MaCompl} and~\cite{PeOrbCompl} we get a spine
of $X_\Sigma$ with at most $c\left(X\right)$ vertices by cutting $P$
open along $\Sigma$.
\end{proof}

\begin{prop}\label{normal:ess:prop}
Let a handle decomposition of a link-pair $X$ be fixed.
If $X$ contains an essential $0$-sphere then it contains a normal one.
If $X$ does not contain essential $0$-spheres but it contains an essential
$2$-sphere then it contains a normal one.
\end{prop}

\begin{proof}
The normalization moves of Haken's theory, besides isotopy relative to $L(X)$,
can be described in a unified fashion as a compression followed by the
removal of one of the two spheres resulting from the compression.
In the context of link-pairs one additionally sees that the disc
along which the compression is performed is disjoint from the link.
Now recall that a $0$-sphere is essential if and only if it is non-trivial.
Of course the compression of a non-trivial $0$-sphere cannot give rise
to two trivial $0$-spheres, and the first conclusion follows.
Now consider an essential $2$-sphere $\Sigma$. Since $\Sigma$
is unsplittable, its compression along a disc disjoint from the link
gives a $0$-sphere and a $2$-sphere. By assumption the first one is trivial,
therefore the second one is not. In addition it is itself
unsplittable, otherwise $\Sigma$ would be, and the conclusion follows.
\end{proof}

\paragraph{Complexity of trivial link-pairs and conclusion}
We are now ready to state and prove the main results of this
section, but we first need the following immediate fact:

\begin{lemma}\label{zero:compl:triv:lem}
Both trivial link-pairs have complexity $0$.
\end{lemma}

\begin{thm}\label{root:add:thm}
If $Y$ is any $\left(0,2\right)$-root of a link-pair $X$ then $c\left(Y\right)=c\left(X\right)$.
\end{thm}

\begin{proof}
Proposition~\ref{normal:ess:prop} implies that we can get a root $Z$ of $X$ by
successively doing surgery on $X$ along $i$-spheres that are normal
with respect to spines having as many vertices as the complexity.
Propositions~\ref{superadd:prop} and~\ref{subadd:prop} then imply that
$c\left(Z\right)=c\left(X\right)$. Now by Proposition~\ref{quasi:unique:02:root:prop}
we have that $Z$ and $Y$ differ only for trivial link-pairs, and
the complexity of a disconnected link-pair is of course the
sum of the complexities of the connected components, therefore
the conclusion follows from Lemma~\ref{zero:compl:triv:lem}.
\end{proof}

\begin{thm}\label{0:add:thm}
Given any two link-pairs $X_1,X_2$ and a $0$-connected sum
$X=X_1\#_0X_2$ one has $c\left(X\right)=c\left(X_1\right)+c\left(X_2\right)$.
\end{thm}

\begin{proof}
Let $\Sigma$ be the $0$-sphere in $X$ along which the $\#_0$ was performed.
If $\Sigma$ is inessential (\emph{i.e.}, trivial) then
up to switching indices we have that $X_2$ is trivial, so $c\left(X_2\right)=0$,
and $X=X_1$, so the conclusion follows. Otherwise doing surgery
along $\Sigma$ in the first place we see that $X$ and $X_1\sqcup X_2$ have a common
root, and the conclusion follows from Theorem~\ref{root:add:thm}.
\end{proof}

\begin{thm}\label{2:quasi:add:thm}
Given any two link-pairs $X_1,X_2$ not containing $1$-spheres
and a $2$-connected sum
$X=X_1\#_2X_2$ one has $c\left(X\right)=c\left(X_1\right)+c\left(X_2\right)$.
\end{thm}

\begin{proof}
The idea of the proof is the same as above, except that now the $2$-sphere
$\Sigma$ in $X$ giving the $\#_2$ \emph{a priori} has two ways of being
inessential: either it is splittable or it is trivial. The first case
is however ruled out by the assumption that $X_1$ and $X_2$ do not
contain $1$-spheres, and the conclusion easily follows.
\end{proof}

We will now explain why the restriction in
Theorem~\ref{2:quasi:add:thm} cannot be avoided. From now on we
denote by $D$ the link-pair $\left(S^2\times S^1,\{*\}\times
S^1\right)$. Note that $c\left(D\right)=0$ because $D$ has a spine
of the form $\left(\{*'\}\times S^1\right)\cup\left(S^2\times
\{*\}\right)$.

\begin{lemma}\label{1-sphere:then_D:lem}
If a link-pair $X$ contains a $1$-sphere then there exist a link-pair $Y$ and $k\geqslant 1$
such that $X$ is the $0$-connected sum of $Y$ and $k$ copies of $D$, and $Y$ does not contain $1$-spheres.
\end{lemma}

\begin{proof}
Let $\Sigma$ be a sphere in $X$ meeting once a component $K_1$ of $L\left(X\right)$,
and no other one. The boundary of a regular neighbourhood of $\Sigma\cup K_1$ is a $0$-sphere,
doing surgery along which we get some $Y_1$ and $D$. If in $Y_1$ there still is a $1$-sphere
we proceed, and we stop in a finite number of steps because $H_1\left(M\left(X\right);\matQ\right)$ has finite dimension.
\end{proof}

Now suppose we have a $2$-connected sum $X=X_1\#_2 X_2$. If we express $X_j$ as the $0$-sum of $Y_j$ and
$k_j$-copies of $D$, as in the previous lemma, we have now three cases for the two
components of $L\left(X_1\right)$ and $L\left(X_2\right)$ involved in the $2$-sum:
\begin{itemize}
\item They belong to $Y_1$ and $Y_2$;
\item They belong to two of the $D$'s;
\item One of them belongs to a $Y_j$ and one to a $D$.
\end{itemize}
Theorems~\ref{0:add:thm} and~\ref{2:quasi:add:thm} easily imply that
$c\left(X\right)=c\left(X_1\right)+c\left(X_2\right)$ in the first
case. The same is true in the second case, because by
Proposition~\ref{superadd:prop} we have $c\left(D\#_2
D\right)\leqslant c\left(D\right)+c\left(D\right)=0+0=0$, so
$c\left(D\#_2 D\right)=0$. However complexity can decrease
arbitrarily in the third case, as the next result implies:

\begin{prop}\label{any:knot:plus:D:is:D:prop}
For every knot $K$ in $S^3$ one has $\left(S^3,K\right)\#_2 D=D$.
\end{prop}

\begin{proof}
By Proposition~\ref{1-sphere:then_D:lem} we know that
$\left(S^3,K\right)\#_2 D=Y\#_0D$ for some $Y$. From the construction, one sees that
the complement of a ball in $Y$ is obtained by digging a tunnel along a properly
embedded arc in $S^2\times[0,1]$ with its ends on the two boundary components.
Realizing $S^2\times[0,1]$ as a subset of $S^3$ one readily sees that the complement
is actually a ball, which implies that $Y=\left(S^3,\emptyset\right)$, and the conclusion follows.
\end{proof}

Since there are infinitely many prime knots $K$ in $S^3$,
Theorem~\ref{natural:thm} implies that $c\left(S^3,K\right)$ attains arbitrarily large values,
so the previous result implies that there cannot exist any lower bound on
$c\left(X_1\#_2 X_2\right)$ in terms of $c\left(X_1\right)$ and $c\left(X_2\right)$.

We conclude with a fact that easily follows from the existence and virtual uniqueness of $(0,2)$-roots, from
Lemma~\ref{1-sphere:then_D:lem} and from Theorem~\ref{natural:thm}, and which (in our opinion) adds to the value of
Theorem~\ref{root:add:thm}:

\begin{prop}
Let $Y$ be a root of a link-pair $X$. Then $Y=Z\sqcup W$ where $Z=Z_1\sqcup\ldots\sqcup Z_n$ and $W=W_1\sqcup\ldots\sqcup W_m$, and:
\begin{itemize}
\item $c(X)=\sum_{j=1}^nc\left(Z_j\right)$ and each $c\left(Z_j\right)$ is the minimal number of
tetrahedra needed to triangulate $Z_j$;
\item Each $c\left(W_j\right)$ vanishes and $W_j$ is
either $D=\left(S^2\times S^1,\{*\}\times S^1\right)$ or one of the $6$ pairs described in
Theorem~\ref{natural:thm}.
\end{itemize}
\end{prop}

\section{Complexity \emph{vs.} the crossing number}

In this section we restrict ourselves to link-pairs $X$ with
$M(X)=S^3$, and we abbreviate the notation indicating $X$ by $L:=L(X)$ only.
We will discuss relations between the complexity $c(L)$ of
$L$ and its crossing number $\crn(L)$. We will also denote by
$\#L$ the number of components of $L$.

\paragraph{A linear upper bound on complexity}
Our first result is based on an easy explicit construction:

\begin{prop}\label{upper:link-compl:prop}
For any link $L\subset S^3$ one has
$$c(L)< 4\crn(L)+2\cdot\#L.$$
\end{prop}

\begin{proof}
Given a diagram $D$ of $L$ on a 2-sphere
$S^2\subset S^3$ such that $D$ realizes the crossing number
of $L$, we can construct a spine of $L$ as follows.
We first dig a tunnel in $S^2$ along each component of $L$,
as suggested in Fig.~\ref{crossing-point:fig} near a crossing of $D$.
    \begin{figure}
    \begin{center}
    \includegraphics[scale=0.5]{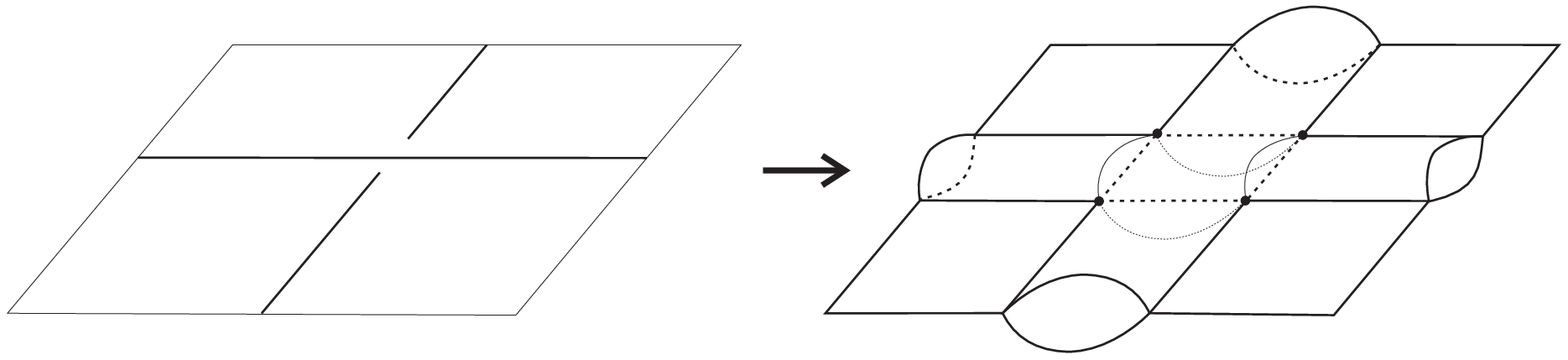}
    \caption{Digging tunnels in $S^2$ along a link diagram}
    \label{crossing-point:fig}
    \end{center}
    \end{figure}
This gives a quasi-special polyhedron
$Q\subset S^3$ whose complement in $S^3$ consists of two balls and a
regular neighbourhood of $L$. For each component of $L$ we then add
to $P$ a 2-disc as in Fig.~\ref{transverse-disc:fig}, getting a
spine $P$
    \begin{figure}
    \begin{center}
    \includegraphics[scale=0.5]{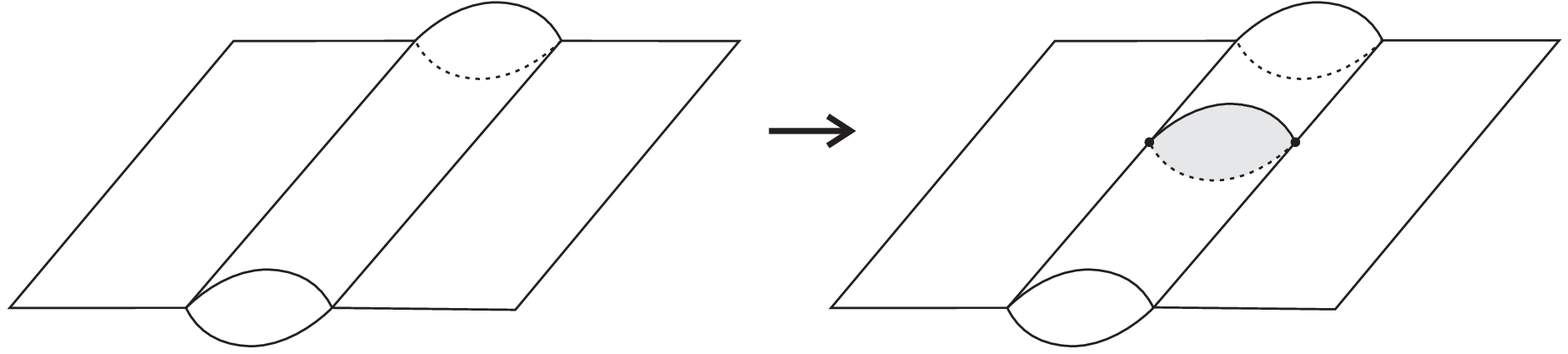}
    \caption{Gluing a transverse disc}
    \label{transverse-disc:fig}
    \end{center}
    \end{figure}
of $L$ with $4\crn(L)+2\cdot\#L$ vertices.
Moreover $P$ can be punctured in two points and collapsed,
so $c(L)$ is strictly smaller than $4\crn(L)+2\cdot\#L$.
\end{proof}

\paragraph{Lower bounds}
Turning to lower bounds, we first prove
the following:

\begin{prop}\label{lower-branched:prop}
If $L$ is a link in $S^3$ and $V$ is
the double cover of $S^3$ branched along $L$, then
$c(L)\geqslant\frac12 c(V)$.
\end{prop}

\begin{proof}
Let $P$ be a minimal spine of $L$.
Without loss of generality we can assume that $P\cap L$ consists of
surface points of $P$. Since the double cover of a disk branched at one point is again
a disc, and the double cover of a 3-dimensional ball branched along
one unknotted arc is again a ball, lifting $P$ to $V$
we get a spine $Q$ of $V$ with twice as many vertices as $P$.
\end{proof}

\begin{cor}\label{homology:bound:cor}
If $L\subset S^3$ is a prime non-split link and
$V$ is the double cover of $S^3$ branched along $L$, then
$$c(L)>\log_5\left|\Tor(H_1(V))\right|-1.$$
\end{cor}

\begin{proof}
The assumptions imply that $L$ is $(0,1,2)$-irreducible in the sense
of link-pairs. The conclusion then follows from
Proposition~\ref{lower-branched:prop}, the irreducibility of $V$
established in~\cite[Corollary 4]{KimToll}, and the lower bound on
the complexity of irreducible 3-manifolds proved
in~\cite[Theorem~1]{MaPe}; note that if $V$ is one of $S^3$,
$S^2\times S^1$, $\matP^3$, $L(3,1)$, which are not covered by~\cite[Theorem~1]{MaPe},
the right-hand side of the inequality in the statement attains negative values,
so the inequality holds in these cases too.
\end{proof}

\begin{thm}\label{Alex:bound:knots:thm}
If $K\subset S^3$ is a prime knot, then
$$c(K)>\log_5\left|\Delta_K(-1)\right|-1,$$
where $\Delta_K(t)$ is the Alexander polynomial of $K$.
\end{thm}

\begin{proof}
Let $V$ be defined as above. Then one knows that $H_1(V)$ is finite,
so it coincides with its torsion part, and
$\left|H_1(V)\right|=\left|\Delta_K(-1)\right|$, as shown for
instance in~\cite[Corollary~9.2]{Lick}. The conclusion then follows
from Corollary~\ref{homology:bound:cor}.
\end{proof}

\paragraph{Two-sided bounds}
Combining the previous results we get:

\begin{thm}\label{c:vs:crn:summary}
If $K\subset S^3$ is a prime alternating knot, then
$$\log_5\left(\crn(K)\right)-1<c(K)< 4\crn(K)+2.$$
\end{thm}

\begin{proof}
The upper bound is the content of
Proposition~\ref{upper:link-compl:prop}. The lower bound follows
from Theorem~\ref{Alex:bound:knots:thm} and the formula
$\crn(K)\leqslant \left|\Delta_K(-1)\right|$, valid for
any alternating knot $K$, as shown for instance in~\cite[Proposition~13.30]{BurdeZies}.
\end{proof}

The next results describe two very specific infinite families of
links for which we can provide more accurate estimates than those
given by Theorem~\ref{c:vs:crn:summary}.
We will show that complexity is asymptotically equivalent to the
logarithm of the crossing number for the links in the first family,
and to the crossing number itself for those in the second family.
Note however that the first family consists of prime but
non-alternating knots, so a direct comparison with
Theorem~\ref{c:vs:crn:summary} is impossible. The second family
consists of prime non-split links, but it contains an infinite
subfamily of knots.

\medskip

To define our first family we consider the Fibonacci numbers
$\left(f_n\right)_{n=0}^\infty$, with $f_0=f_1=1$, and we
denote by $T_{m,q}$ the $(m,q)$ torus knot.

\begin{thm}
If $n\geqslant 4$ and $n\equiv 0,2\pmod 3$, then
$$\frac12\log_5\left( \crn\left(T_{f_n,f_{n-1}}\right)\right)-2<c\left(T_{f_n,f_{n-1}}\right)
<4\log_2\left(\crn\left(T_{f_n,f_{n-1}}\right)\right)-1.$$
\end{thm}

\begin{proof}
We first recall that
$\crn\left(T_{m,q}\right)=\min\{m(q-1),q(m-1)\}$,
as shown for instance in~\cite{TorusCrossing}. It
follows that
$\crn\left(T_{f_n,f_{n-1}}\right)=f_n\cdot\left(f_{n-1}-1\right)$,
hence
$$\log_a f_n<\log_a \left(\crn\left(T_{f_n,f_{n-1}}\right)\right)<2\log_a f_n,$$
and finally
$$\frac12\log_a \left(\crn\left(T_{f_n,f_{n-1}}\right)\right)
<\log_a f_n<\log_a \left(\crn\left(T_{f_n,f_{n-1}}\right)\right).$$
The conclusion is then a consequence of the following claims:

\medskip

\noindent\textsc{Claim:} $c\left(T_{f_n,f_{n-1}}\right)\leqslant
4\log_2f_n-1$. It was shown in~\cite[Section~6.3]{HHMP}
that $c\left(T_{m,q}\right)\leqslant\lambda\left(m,q,0,1\right)-3$ for a certain function
$\lambda$ defined there, and it follows immediately from its
definition that $\lambda\left(m,q,0,1\right)$ is not greater
than twice the sum of the partial quotients in the expansion of
$m/q$ as a continued fraction. This sum is equal to $n$
for $f_n/f_{n-1}$, therefore $c\left(T_{f_n,f_{n-1}}\right)\leqslant 2n-3$.
Moreover one easily sees that
$n\leqslant 2\log_2f_n+1$, whence the claimed inequality.

\medskip

\noindent\textsc{Claim:} $c\left(T_{f_n,f_{n-1}}\right)> \log_5f_n-2$.
Recall that in general
$$\Delta_{T_{m,q}}(t)=\frac{(1-t)(1-t^{mq})}{(1-t^m)(1-t^q)}.$$
If $m$ is even and $q$ is odd we then have
$$\Delta_{T_{m,q}}(-1)=\frac{1-t}{1-t^q}\left(1+t^m+\ldots+t^{m(q-1)}\right)\Big|_{t=-1}=q.$$
Likewise, if $q$ is even and $m$ is odd then $\Delta_{K_{m,q}}(-1)=m$.
Now the assumption that $n\equiv 0\pmod 3$ or $n\equiv 2\pmod 3$ implies that one of
$f_n$ and $f_{n-1}$ is even and the other one is odd, whence
$\Delta_{T_{f_n,f_{n-1}}}(-1)\geqslant f_{n-1}>f_n/2$. The claimed
inequality now easily follows from Theorem~\ref{Alex:bound:knots:thm}.
\end{proof}

Our second family is that of the so-called Turk's head links ${\rm
Th}_n$ (see, for instance~\cite{MedVes}), defined as the closure of
the 3-string braid $\left(\sigma_1\sigma_2^{-1}\right)^n$. Notice
that ${\rm Th}_2$ is the figure-eight knot, ${\rm Th}_3$ is the
Borromean rings, ${\rm Th}_4$ is the Turk's head knot $8_{18}$, and
${\rm Th}_5$ is the knot $10_{123}$.

\begin{thm}
For sufficiently large $n$ one has
$$\frac12 \crn\left({\rm Th}_n\right)\leqslant c\left({\rm Th}_n\right)< 4\crn\left({\rm Th}_n\right)+6.$$
\end{thm}

\begin{proof}
The diagram of ${\rm Th}_n$ obtained by closing
$\left(\sigma_1\sigma_2^{-1}\right)^n$ is reduced alternating and it
has $2n$ crossings, therefore $\crn\left({\rm Th}_n\right)=2n$
by~\cite[Corollary~5.10]{Lick}). The upper estimate on complexity
then follows from Proposition~\ref{upper:link-compl:prop} and the
fact that ${\rm Th}_n$ has at most 3 components.

Now let $V_n$ be the double cover of $S^3$ branched along ${\rm
Th}_n$. By~\cite[Corollary~3.4]{MaPeVe} we have that $V_n$ is the
so-called $n$-th Fibonacci manifold~\cite{FibonacciMfd}. It was
shown in~\cite[Theorem~1]{MedVes} that $c(V_n)\geqslant 2n$ for
sufficiently large $n$, and the conclusion follows from
Proposition~\ref{lower-branched:prop}.
\end{proof}

\section{Complexity \emph{vs.} hyperbolic volume}

In this section we compare the complexity of a link with the
hyperbolic volume (if any) of its exterior. We consider again the
general case of a link-pair $X$, and we define $E(X):=M(X)\setminus
L(X)$, but rather soon we will return to $M(X)=S^3$, in which case
we write $E(L)$ instead of $E(X)$, with $L:=L(X)$.

\begin{prop}\label{lower:volume:prop}
If $X$ is a link-pair with $L(X)\neq\emptyset$ and hyperbolic $E(X)$, then
$$c(X)>\frac1{v_3}{\vol\left(E(X)\right)},$$
where $v_3=1.01494\ldots$ is the volume of the regular ideal tetrahedron.
\end{prop}

\begin{proof}
Being hyperbolic, $E(X)$ is irreducible, boundary-incompressible,
and acylindrical, which easily implies that $X$ is
$(0,1,2)$-irreducible. Therefore a minimal spine $P$ of $X$ is
special by Theorem~\ref{natural:thm}. Puncturing $P$ at its
intersections with $L(X)$ and collapsing we get a spine of $E(X)$
with strictly fewer vertices than $P$, so $c(X)>c\left(E(X)\right)$.
We conclude using the general inequality $\vol(M)\leqslant v_3\cdot
c(M)$, valid for any finite-volume hyperbolic $M$, as already
remarked by Thurston~\cite[Corollary~6.1.7]{bible} (see also~\cite{Anis},
and~\cite[Proposition~2.7]{MaPeVe} for a formalization in the closed
case).
\end{proof}

We now turn to the case of links in $S^3$ and
recall some terminology introduced by Lackenby
in~\cite{La}. A \emph{twist} in a link diagram $D\subset\matR^2$ is
either a maximal collection of bigonal regions of $\matR^2\setminus D$
arranged in a row, or a single crossing with no incident bigonal
regions. The \emph{twist number} $t(D)$ of $D$ is the total number
of twists in $D$. The remarkable achievement of~\cite{La}
was to show that for an alternating link $L\subset S^3$
the volume of $E(L)$ is bounded from above and from below
by degree-1 polynomials (with positive leading
coefficients) in the twist number of the alternating diagram of $L$.

It is now quite easy to construct infinite series of alternating
hyperbolic knots and links for which the twist number differs by a
constant from a non-zero multiple of the crossing number. Combining
Propositions~\ref{upper:link-compl:prop} and~\ref{lower:volume:prop}
with the cited result of~\cite{La} we then get that for any such
series the complexity can be efficiently estimated in terms of the
volume of the exterior. The next result describes an example:

    \begin{figure}
    \begin{center}
    \includegraphics[scale=0.6]{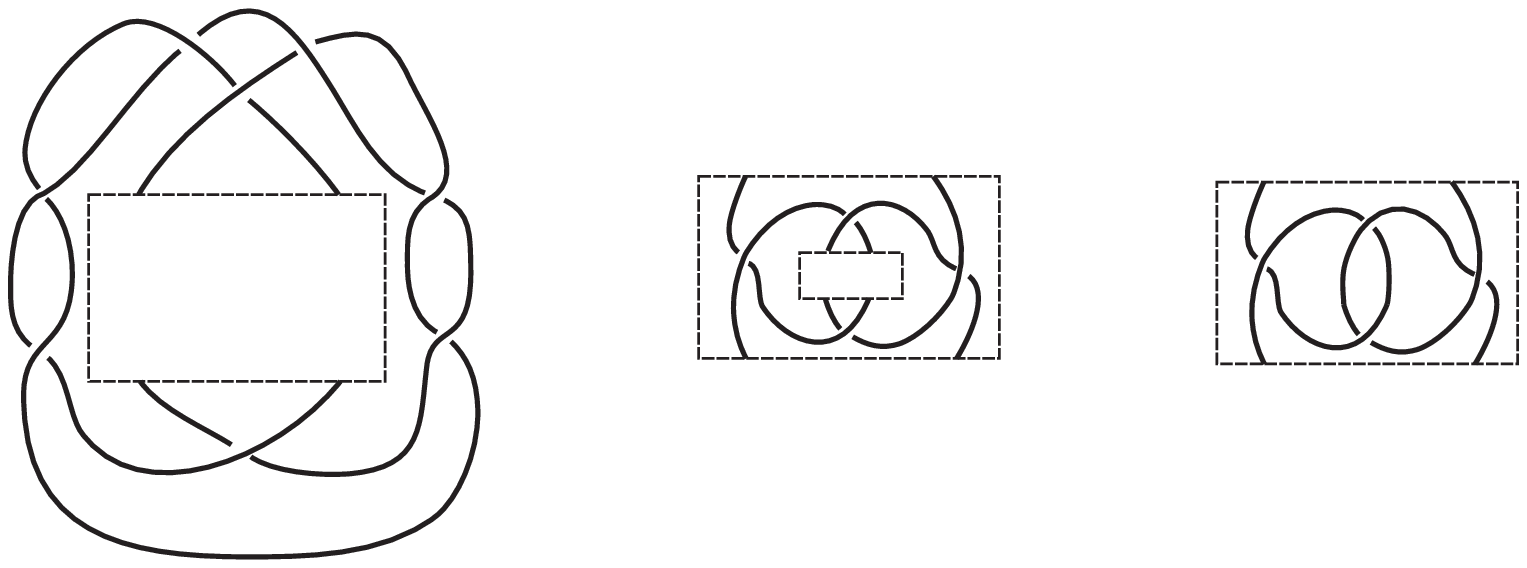}
    \mycap{The knot $K_n$ is obtained starting from the diagram on the left,
    successively inserting in it $n$ copies of the diagram in the centre, and
    concluding with the diagram on the right}
    \label{my-knots:fig}
    \end{center}
    \end{figure}

\begin{prop}\label{2-sided:hyp:prop}
The knots $\left(K_n\right)_{n=0}^{\infty}$ defined in Fig.~\ref{my-knots:fig}
are hyperbolic. Moreover $c\left(K_n\right)$ is bounded from above and from
below by degree-$1$ polynomials with positive leading coefficients of each
and every of the following:
\begin{itemize}
\item $\vol\left(E\left(K_n\right)\right)$;
\item $\crn\left(K_n\right)$;
\item $t\left(K_n\right)$;
\item $n$.
\end{itemize}
\end{prop}

\begin{figure}
    \begin{center}
    \includegraphics[scale=0.6]{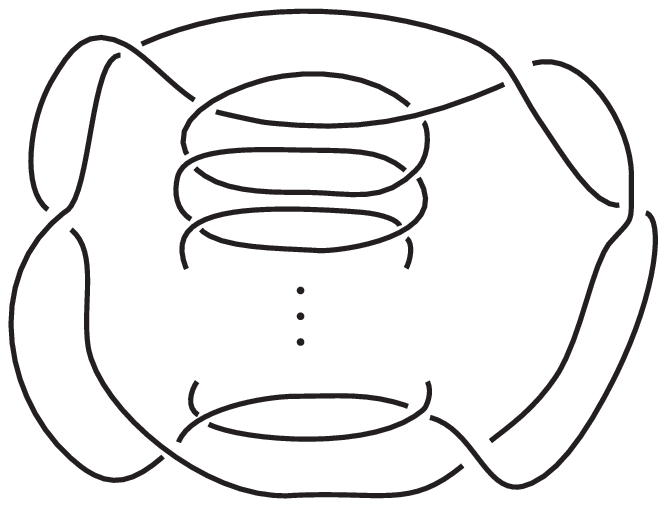}
    \mycap{The link $L_n$, with $n$ standing for the number of central horizontal bigons}
    \label{my-links:fig}
    \end{center}
    \end{figure}

\begin{figure}
    \begin{center}
    \includegraphics[scale=0.5]{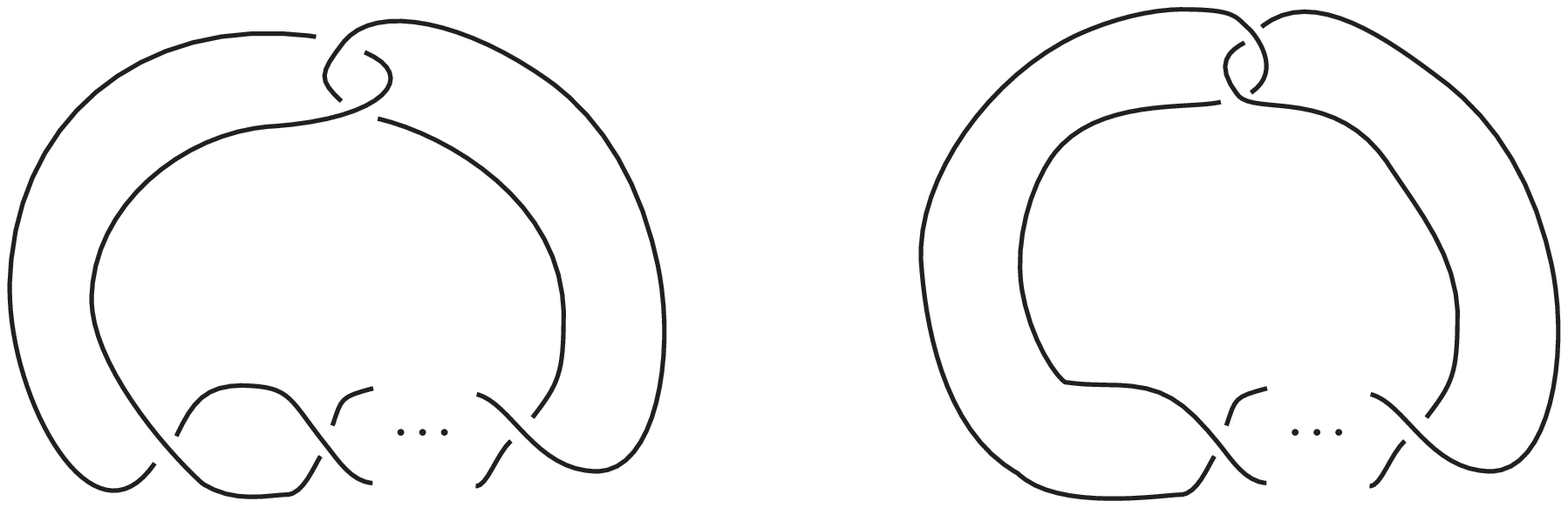}
    \mycap{Left: the usual diagram $n$-th twist knot, with $2n+1$ crossings.
    Right: an alternating diagram of the same knot, showing that its crossing number is $2n$}
    \label{twist-knots:fig}
    \end{center}
    \end{figure}

A similar statement holds for the links $\left(L_n\right)_{n=0}^\infty$
defined in Fig.~\ref{my-links:fig}.
We conclude by noting that however no such result
can hold in general. Considering the twist knots described in Fig.~\ref{twist-knots:fig}
we see that they are hyperbolic and distinct, so their complexity attains
arbitrarily large values. However their volume stays bounded, as one can
easily see using~\cite{La} or hyperbolic Dehn surgery~\cite{BePe}.

\vspace{1cm}

\noindent
Dipartimento di Matematica Applicata\\
Universit\`a di Pisa\\
Via Filippo Buonarroti, 1C\\
56127 PISA -- Italy\\
\ \\
pervova@guest.dma.unipi.it\\
petronio@dm.unipi.it

\end{document}